\documentclass[final]{article}
\usepackage{amssymb,amsmath,amsfonts,amsthm,mathrsfs,graphicx,color}
\usepackage[colorlinks=true,linkcolor=blue]{hyperref}

\usepackage{amsmath,amssymb,amsfonts,mathrsfs,color}
\usepackage{graphicx,epsf}
\newtheorem{theor}{Theorem}
\newtheorem{coro}{Corollary}[section]
\newtheorem{lem}{Lemma}[section]

\newtheorem{remark}{Remark}

\newcommand{\N}{{\mathbb N}}
\newcommand{\R}{{\mathbb R}}
\newcommand{\Z}{{\mathbb Z}}

\newcommand{\calL}{\mathcal L}
\newcommand{\calO}{\mathcal O}
\newcommand{\calV}{\mathcal V}

\newcommand{\dive}{{\rm div}}

\newcommand{\epsl}{\varepsilon}

\newcommand{\refq}[1]{~(\ref{#1})}
\newcommand{\rmd}{\mbox{d}}
\newcommand{\sauf}{\setminus}

\newcommand{\teta}{\theta}
\newcommand{\Om}{\Omega}

\newcommand{\ove}{\overline}

\title
      {Recovering time-dependent inclusion in heat conductive bodies by a dynamical probe method}
\title{Recovering time-dependent inclusion in heat conductive bodies by a dynamical probe
 method}

\author{O. Poisson\thanks{Aix Marseille Universit{\'e}, I2M, UMR CNRS 6632,
 France ({\tt  olivier.poisson@univ-amu.fr}).}}

\begin{document}
\maketitle

%

\begin{abstract}
 We consider an inverse boundary value problem for the heat equation $\partial_t v =
 \dive_x\,(\gamma\nabla_x v)$ in $(0,T)\times\Om$, where $\Om$ is a bounded  domain of $\R^3$,
 the heat conductivity $\gamma(t,x)$ admits a surface of discontinuity which depends on time
 and without any spatial smoothness.

%
 The reconstruction and, implicitly, uniqueness of the moving inclusion, from the knowledge
 of the Dirichlet-to-Neumann operator,  is realised by a dynamical probe method based on
 the construction of fundamental solutions of the elliptic operator $-\Delta + \tau^2\cdot$,
 where $\tau$ is a large real parameter, and a couple of inequalities relating data
 and integrals on the inclusion, which are similar to the elliptic case.
 That these solutions depend not only on the pole of the fundamental solution,
 but on the large parameter $\tau$ also, allows the method to work in the very general situation.
\end{abstract}

\medskip
 {\bf Keywords}: Inverse problem, Dirichlet-to-Neumann map, heat probing.

\medskip
{\bf AMS} : 35R30, 35K05.

\pagestyle{myheadings}
\thispagestyle{plain}


\section{Introduction}
\subsection{Inverse heat conductivity problem}
 Let $\Om$ be a bounded domain in $\R^3$, with lipschitzian boundary $\Gamma = \partial\Om$,
 and consider the following initial boundary value problem
\begin{eqnarray}
\left\{ \begin{array}{rll}
 \partial_t v &=& \dive_x\, (\gamma \nabla_x v)  \quad {\rm in} \quad \Om_T= (0,T)\times\Om,\\
 v &=& f  \quad {\rm on} \quad \Gamma_T=(0,T)\times\Gamma, \\
 v\big|_{t=0} &=& v_0  \quad {\rm on} \quad \Om,
\end{array} \right.
\label{parab.v}
\end{eqnarray}
 where $\gamma = \gamma(t,x)\in W^{1,\infty}((0,T);L^\infty(\Om))$ with the following properties:

\medskip
(C-$\gamma$) {\it There exist a positive function $(t,x)\mapsto k(t,x)$ and, for all $t\in [0,T]$,
 a non empty open set $D(t) \subset \Om$}, such that
\begin{itemize}
\item $\frac1C\le k \le C\quad {\rm in } \quad D_T:=\cup_{[0,T]}\{t\}\times D(t)$ for some $C>1$,
\item $k-1$ has a constant sign in $D_T$,
\item
\begin{eqnarray*}
 \gamma(t,x) = \left\{ 
\begin{array}{ll}
 1 \quad {\rm if} \quad x \not\in D(t),\\
 k(t,x) \quad {\rm if} \quad  x\in D(t).
\end{array} \right.
\end{eqnarray*}
\end{itemize}
 {\em We don't assume any smoothness neither on $D(t)$ nor  that
  $\partial D(t)\cap\Gamma=\emptyset$}.

\medskip
 Our main purpose is to study discontinuous perturbations, however, we allow
 $\gamma(t,\cdot)$ to be continuous. Hence, we impose the following assumption.

\medskip
\noindent
(C-D) \ \  
 $\inf_{x\in K} |k(t,x)-1| > 0$  {\it for any compact set} $K \subset D(t)$,
 for all $t\in[0,T]$.
 
\medskip 
 We shall consider a large parameter $\tau>0$ and allow the initial data $v_0(x)$
 to depend on $\tau$, under the following condition:

\medskip
\noindent
(C-0) \ \
 There exist $\tau$-independent positive constants $C$, $l_0$ such that
 
 $\|v_0\|_{L^2(\Om)} \le Ce^{\tau l_0}$, for all $\tau$.

\medskip
\noindent

\bigskip
 Physically, the region $D(t)$  corresponds to some inclusion in the medium with heat
 conductivity different from that in the background  domain $\Om$.
 The  problem we address in this paper is to determine  $D(t)$ by using the  knowledge of
 the Dirichlet-to-Neumann map (D-N map) :
$$
 \Lambda_{\gamma,v_0} : f \mapsto \partial_{\nu} v(t,x), \quad (t,x) \in \Gamma_T,
$$
 where $v=\calV(\gamma;f)$ denotes the unique solution of (\ref{parab.v}),
 $\nu$ is the outer unit normal to $\Gamma$, and
 $\partial_\nu = \dfrac{\partial}{\partial\nu} = \nu\cdot\nabla_x$.
 In physical terms, $f=f(t,x)$ is the temperature distribution on the boundary and
 $\Lambda_{\gamma,v_0}(f)$ is the resulting heat flux through the boundary.
 
 The above inverse boundary value problem is related to nondestructive testing where one looks
 for anomalous materials inside a known material.

\medskip
 To clarify our purpose, we remind briefly Ikehata's probe method for the elliptic inverse problem.
\subsection{The elliptic situation}
 In the probe method for the well-known elliptic situation, Problem\refq{parab.v} is replaced by
\begin{eqnarray}
\label{ellip.v}
\left\{ \begin{array}{rll}
 \dive_x\, (\gamma \nabla_x v) &=& 0 \quad {\rm in} \quad \Om,\\
 v &=& f  \quad {\rm on} \quad \Gamma, \\
\end{array} \right.
\end{eqnarray}
 $D_T$ is replaced by an open set $D\subset\Om$.
 The Dirichlet-to-Neumann operator $\Lambda_{\gamma}$ is a mapping:
 $H^{\frac12}(\Gamma)\ni f \mapsto \partial_{\nu} v\in H^{-\frac12}(\Gamma)$, where $v$ is
 the unique solution of\refq{ellip.v}.
 The probe method (see \cite{IKE.PRO}) starts by considering the fundamental solution $h_0(x)= \frac1{4\pi |x-y|}$
 of $-\Delta h_0=\delta_y$, with pole $y\in\Om$ .
  Then, one approximates $h_0$ outside a needle $\Sigma\subset \ove\Om$ with one end on $\Gamma$
  and the other one being $y$ by a sequence $\{h_j\}_{j\ge 1}$ such that $-\Delta_x h_j=0$ in $\Om$,
  and estimates $\int_D |\nabla h_j(x)|^2 dx$ (or $\int_D |\nabla h_0(x)|^2 dx$) thanks to the following
  couple of inequalities:
\begin{eqnarray}
\label{i1.ell}
 \frac1C \int_D |\nabla h_j(x)|^2 \le 
 \left|\int_{\Gamma} ( \Lambda_{\gamma} (h_j|_{\Gamma})-\partial_\nu h_j)
  h_j|_{\Gamma} \right| \le  C \int_{D} |\nabla h_j(x)|^2,
\end{eqnarray}
 for some $C>1$ which does not depend on $h_j$.
 However, in the parabolic situation, inequalities like\refq{i1.ell} are unclear, except
 in the static case $D(t)=D(0)$, for a reduced class of functions as $h(t,x)=e^{\tau^2 t}p(x)$,
 where $\tau>0$ is a large parameter and $p$ satisfies $-\Delta p +\tau^2 p = 0$ in $\Om$.
 See \cite{GAI.PRO}, and \cite{IKE.ENC,IKE.REC} also for a similar approach.

\medskip
 In our work we built fundamental solutions for parabolic operators which slightly differs from
 the heat operator $\calL_1 \equiv \partial_t - \Delta_x$ or to its adjoint
 $\calL_1^* \equiv -\partial_t - \Delta_x$.
 Then, we obtain a couple of  inequalities corresponding to\refq{i1.ell}, and,
 thanks to that, we run the probe method for the reconstruction of $D_T$.
\subsection{Special fundamental solutions}
 Let $\tau>0$, $y\in\R^3$, we consider the following function
$$
 p_{\tau,y}(x) := \frac1{r} e^{-\tau r},\quad r=|x-y|.
$$
 It satisfies $p_{\tau,y}\in L^2(\R^3)$ and
\begin{equation}
\label{e.ptau}
  -\Delta p_\tau +\tau^2 p_\tau = \delta(x-y)  \mbox { in $\R^3$.}
\end{equation}
 Thus, the functions $h_\tau(t,x)=e^{\tau^2t}p_\tau(t,x)$ and  
 $h^*_\tau(t,x)=e^{-\tau^2t}p_\tau(t,x)$ are respectively solutions to the following
 parabolic equations
\begin{equation}
\label{e.htau}
 \calL_1 h := \partial_t h- \Delta_x h = \delta(x-y), \quad
 \calL_1^* h^* := -\partial_t h^*- \Delta_x h^* =  \delta(x-y).
\end{equation}
 Of course, since the right term in Equalities\refq{e.htau} is not $\delta(t-s)\otimes\delta(x-y)$,
 our function $h_\tau$ differs from the usual fundamental solution
 $G_{s,y}(t,x) = \frac{\chi_{(s,+\infty)}(t)}{(4\pi (t-s))^{3/2}} e^{-\frac{|x-y|^2}{4(t-s)}}$,
 with parameter $(s,y)\in\R\times\R^3$, considered in \cite{CRI.STA,DAI.PRO,ISA.REC}.

\medskip
 Let $\Om'$ be any smooth bounded domain containing $\ove\Om$.
 Let $\Sigma$ be an finite or infinite globally lipschitzian curve in $\Om'$
 with parameter $y=y(t)$, $t\in\R$, such that $y(t)\not\in\ove\Om$ for $t\le 0$.
 We shall use the restriction $\Sigma|_{-1\le t \le T+1}$ only, which plays the role
 of a Ikehata's needle, but starts at point $y(-1)$ outside $\ove\Om$ and can be self-intersecting.
 We put
$$ p_{\tau,y(t)} =: p_{\tau,\Sigma} (t,x),\quad (t,x)\in \R\times\R^3,$$ 
 and replace the functions $h_\tau$,$ h_\tau^*$, respectively by solutions of
\begin{eqnarray}
\label{p.U}
 (\calL_1 +\tau\eta(t)) U(t,x) &=& \delta(x-y(t))  \quad {\rm in} \quad \R\times\Om,\\
\label{p.U*}
 (\calL_1^*+\tau\eta(t)) U^*(t,x) &=& \delta(x-y(t))  \quad {\rm in} \quad \R\times\Om,
\end{eqnarray}
 where functions $\eta(t)$ are not necessarily the same in\refq{p.U} and\refq{p.U*},
 but belong to $L^\infty(\R)$ with a $\tau$- independent upper bound:
\begin{equation} 
\label{est.eta}
 |\eta(t)|  \le \mu,\quad t\in \R,
\end{equation}
%
 for some $\mu>0$ that we shall precise later.
 (By $\dot q$ we denote the derivative of a function $q$ according to the time variable $t$).
  
 The solutions $U$ and $U^*$ we look for should be respectively written $U_{\tau,T,\Sigma}$,
 $U_{\tau,T,\Sigma}^*$, but we simply denote them by $U_\tau$, $U_\tau^*$.
 We construct them so that they have the following form:
%
$$
 U_\tau(t,x) = e^{\tau^2 t} u_\tau(t,x), \quad U_\tau^*(t,x) =
   e^{-\tau^2 t} u_\tau^*(t,x) \quad {\rm in} \quad \Om_T,
$$
%
 with functions $u_\tau(t,x)$, $u^*_\tau(t,x)$ sufficiently close to
 $p_{\tau,\Sigma}(t,x)$ in the sense below.

%

\medskip
\begin{lem}
\label{l.varphi}
 There exists a function $\varphi(t,x;\tau) \in C([0,T];H^1_{loc}(\R^3))$ satisfying
 the following points.
%

%
\begin{equation}
\label{v.varphi0}
 \varphi(t,y(t);\tau)=1,\quad t\in\R,\quad \tau\ge 0.
\end{equation}
%
%
 There exists $\tau_1(\mu)$ such that for all $\tau>\tau_1$, we have
\begin{eqnarray}
\label{est.varphi}
 \quad \frac1{C(\mu)}\le \varphi(t,x;\tau) \le C(\mu) \quad (t,x)\in \R\times \Om,\\
\label{est.dvarphi}
 |\nabla_x\varphi(t,x;\tau)|  \le C(\mu), \quad (t,x)\in \R\times \Om, 
\end{eqnarray}
 for some $C(\mu)>0$ which does not depend on $\tau$ or $(t,x)\in \R\times\Om$.
%
 The function
\begin{equation}
\label{d.Utau}
 U_\tau(t,x):=e^{\tau^2 t}u_\tau(t,x)\equiv  e^{\tau^2 t} \varphi(t,x;\tau) p_{\tau,\Sigma}(t,x)
\end{equation}
 satisfies\refq{p.U} for all $\tau> 0$.
%
\end{lem} 
\medskip
%
\begin{remark}
\label{r.p.U*}
 Similarly we construct $\varphi^*(t,x;\tau) \in C([0,T];H^1_{loc}(\R^3))$
 such that the function
%
$$
 U_\tau^*(t,x) = e^{-\tau^2 t} u_\tau^*(t,x) = e^{-\tau^2 t}  \varphi^*(t,x;\tau)
   p_{\tau,\Sigma}(t,x)
$$
%
 satisfies\refq{p.U*}, with Relation \ref{v.varphi0} and Estimates\refq{est.varphi},
\refq{est.dvarphi} where $\varphi$ is replaced by  $\varphi^*$.
\end{remark}

\medskip
 We choose functions $\eta$ as follows. Let us consider the following function
\begin{equation}
\label{d.kappa}
 \kappa(t) := e^{-\tau\mu |t-\teta|},\quad t\in\R,
\end{equation}
 where the parameters $\mu>0$, $\teta\in [0,T]$, are independent of $\tau$ and will be
 precised later. For $U_\tau$, we put simply $\eta\equiv 0$. For $U_\tau^*$, we put
%
$$
 \eta(t) =-\frac{\dot\kappa}{\tau\kappa} = \mu\, {\rm sgn}(t-\teta).
$$
%
 (Thus each $\eta$ satisfies\refq{est.eta} obviously).

\subsection{Needle sequences}
\label{s.NS0}
 Let $T'\in(0,T]$, we put $\Sigma_{T'}= \{(t,y(t))$, $0\le t\le T'\}$.
 We consider needle sequences $\{U_j\}_j$, $\{U_{j,T'}^*\}_j$, respectively associated
 to $(\Sigma_T, U_\tau, \eta\equiv 0)$, $(\Sigma_{T'}, U^*_\tau,\eta\equiv -\frac{\dot\kappa}{\tau\kappa})$.
 They satisfy:
\begin{eqnarray*}
 (\calL_1+ \tau\eta(t)) U_j & = & 0 \quad{\rm in}\quad \Om_T,\\
 U_j(0,\cdot)  & = & U_\tau(0,\cdot) \quad{\rm in} \quad\Om,\hspace{-1cm}
\end{eqnarray*}
\begin{eqnarray}
\label{e.Uj*}
 (\calL_1^*+ \tau\eta(t)) U_ {j,T'}^* & = & 0 \quad{\rm in}\quad \Om_{T'},\\
\label{ci.Uj*}
 U_{j,T'}^*(T',\cdot)  & = &  U_\tau^*(T',\cdot) \quad{\rm in} \quad \Om,\hspace{-1cm}
\end{eqnarray}
 and, for all open set $V\subset\Om_{T'}$ such that d$(V,\Sigma_{T'})>0$,
\begin{eqnarray}
\label{lim.UjUj*}
 (U_j,U_{j,T'}^*)  &  \stackrel{ j\to \infty}\to &  (U_\tau,U_{\tau}^*) \quad
 \mbox{(strongly) in}\quad (H^{1,0}(V)\cap H^{0,1}(V))^2.
\end{eqnarray}
 Here, we use the notation
$$ H^{s,p}(V) = \ove{C_0^\infty(V)}^{_{\|\cdot\|_{s,p}}}, \quad s,p\ge 0,$$
 with $ \|u\|_{s,p}\equiv \|(t\mapsto \|u(t,\cdot)\|_{H^p(\R^3;dx)})\|_{H^s(\R;dt)}$.
%

\medskip
\begin{remark}
\label{rem.Uj*}
 Since $y(0)\not\in\ove{\Om}$, thanks to\refq{lim.UjUj*}, we have
$$
  U_{j,T'}^*(0,\cdot)  \stackrel{ j\to \infty}\to  U_\tau^*(0,\cdot) \quad \mbox{(strongly) in}
 \quad L^2(\Om). 
$$
\end{remark}

\subsection{Reflected waves}
 Let $V_j:=\calV(\gamma,U_j|_{\Gamma_T})$ be the solution of\refq{parab.v} with data
 $f=U_j|_{\Gamma_T}$, put $W_j=V_j-U_j$.
 By denoting $\calL_\gamma = \partial_t -\nabla_x((\gamma-1)\nabla_x$, we have
\begin{eqnarray}
\label{e.Wj}
 \calL_\gamma W_j &=& \nabla_x((\gamma-1)\nabla_x U_j) \quad{\rm in} \quad \Om_T,\\
\label{bc.Wj}
 W_j &=& 0 \quad{\rm on}\quad \Gamma_T,\\
\nonumber
 W_j|_{t=0} &=& v_0 - U_\tau  \quad{\rm in}\quad \Om.
\end{eqnarray}
 If $\Sigma_T\cap \ove{D_T}=\emptyset$ then, in view of\refq{e.Wj}, we have
\begin{eqnarray}
\label{lim.Wj}
 W_j  \stackrel{ j\to \infty}\to W_\tau \quad {\rm in} \quad C([0,T],H^1(\Om)),
\end{eqnarray}
 where $W_\tau$ is the unique solution of
\begin{eqnarray}
\nonumber
 \calL_\gamma W_\tau &=& \nabla_x((\gamma-1)\nabla_x U_\tau) \quad{\rm in} \quad \Om_T,\\
\label{bc.Wtau}
 W_\tau &=& 0 \quad{\rm on}\quad \Gamma_T,\\
\nonumber
 W_\tau|_{t=0} &=& v_0 - U_\tau  \quad{\rm in}\quad \Om.
\end{eqnarray}
 Thus
\begin{equation}
\label{lim.Vj}
 V_j  \stackrel{ j\to \infty}\to V_\tau = U_\tau +W_\tau\quad {\rm in} \quad H^{0,1}(V),
\end{equation}
 for all open set $V$ such that $\ove V\subset \ove{\Om_T}\sauf\Sigma_T$.
 The function $V_\tau$ satisfies
\begin{eqnarray*}
\nonumber
 \calL_\gamma V_\tau &=& \delta(x-y(t)) \quad{\rm in} \quad \Om_T,\\
 V_\tau &=& U_\tau \quad{\rm on}\quad \Gamma_T,\\
 V_\tau|_{t=0} &=& v_0   \quad{\rm in}\quad \Om.
\end{eqnarray*}
%


\subsection{Pre-indicator sequence and indicator function}
 Let $\teta\in (0,T)$, $\mu>0$, and put
\begin{eqnarray*}
 I_j(\tau,\mu,\teta,T) &:=& \int_{\Gamma_T} (\partial_\nu W_j) \, U_{j,T}^* \, d\sigma(x)
   \kappa(t) dt,\\
 I_\infty(\tau,\mu,\teta,T) &:=&  \int_{\Om_T}(\gamma-1) \nabla_x V_\tau
  \nabla_x U_\tau^* \,  dx \kappa(t) dt + \int_\Om \left[\kappa W_\tau  U_\tau^*\right]^T_0 dx,
\end{eqnarray*}
 where $\kappa$ is defined by\refq{d.kappa}, and $d\sigma(x)$ is the usual measure
 on the boundary $\Gamma$.

\medskip
\begin{remark}
 As opposite to dynamical probe methods in \cite{DAI.PRO,IKE.REC}, the indicator function
 $I_\infty$ depends on the needle $\Sigma$.
\end{remark}

\medskip
 The following result ensures us that we have the knowledge of $I_\infty(\tau,\mu,\teta,T)$
 from the Cauchy data of\refq{parab.v} when $\Sigma_T$ does not touch $D_T$.

\medskip
\begin{lem}
\label{l.0}
 Assume that $\Sigma_T\cap \ove{D_T}=\emptyset$. Then, for all $\teta\in [0,T]$,
 $\mu,\tau\ge 0$, we have
\begin{equation}
\label{lim1.Ij}
 I_j(\tau,\mu,\teta,T) \stackrel{ j\to \infty}\to I_\infty(\tau,\mu,\teta,T) \in\R.
\end{equation}
\end{lem}

\subsection{Main theorems}
 The following theorem separates the cases  $\Sigma_T\cap \ove{D_T} =\emptyset$ and
 $\Sigma_T\cap \ove{D_T}\neq\emptyset$. In the second case, we can put
$$ T^*:=\min \{T'\in [0,T]; \: \Sigma_{T'}\cap \ove{D_{T'}}\neq\emptyset\} >0. $$
 We put also
$$
 \mu_{T',\teta} :=l_0 \max \left(\frac1{\teta},\frac1{(T'-\teta)}\right),
 \quad 0<\teta<T'\le T,
$$
 where $l_0$ is the constant in (C-0).

\medskip
 By $\rmd(y,X)$ we denote the "distance" from the point $y$ to the set $X$,
 and by $\rmd(Y,X)$ the "distance" from the set $Y$ to the set $X$:
$$  \rmd(Y,X) = \inf_{y\in Y} \rmd(y,X),\quad  \rmd(y,X) = \inf_{x\in X} |x-y|. $$
\medskip
\begin{theor}
\label{t.1}
 Assume (C-D) and let $\teta\in (0,T)$. Assume that $\Sigma_T\cap \ove{D_T}=\emptyset$.\\
 Then there exists $\mu_1>\mu_{T,\teta}$ such that, for all $\mu>\mu_1$, there exist
 $C(\mu)>1$, $\tau_1>\mu_1$ such that, under assumption (C-0), we have
\begin{eqnarray}
\label{dine.Iinfty}
 \frac1{C(\mu)} \int_{\Om_T} \!\!\!\! |\gamma-1|\, |\nabla_x p_{\tau,\Sigma}|^2 dx \kappa dt
  \le  \left| I_\infty(\tau,\mu,\teta,T)\right| \le \\
\nonumber
 C(\mu) \int_{\Om_T}\!\!\!\! |\gamma-1|\, |\nabla_x p_{\tau,\Sigma}|^2 dx \kappa dt,
\end{eqnarray}
\end{theor}

\medskip
\begin{theor}
\label{t.2}
 Assume (C-D). Then the following points hold.
 
\medskip
\noindent
 $\bullet$ Let $\teta\in (0,T)$. For all $\mu>\mu_1$, under assumptions (C-0),
 $\Sigma_T\cap \ove{D_T}=\emptyset$, we have
\begin{equation}
\label{lim1.Iinfty}
 \lim_{\tau\to\infty} \frac1{\tau} \ln |\, I_\infty(\tau,\mu,\teta,T) | =
  -2\rmd(y(\teta),D(\teta)) <0.
\end{equation}
$\bullet$  under assumptions (C-0), $\Sigma_T\cap \ove{D_T}=\emptyset$,
 we have
\begin{equation}
\label{lim2.Iinfty}
 \limsup_{\teta\in(0,T)} \lim_{\mu\to\infty}\lim_{\tau\to\infty} \frac1{\tau}
 \ln| \, I_\infty(\tau,\mu,\teta,T)| = -2\rmd(\Sigma_T,D_T) <0.
\end{equation}
\medskip
$\bullet$  Assume that $\Sigma_T\cap \ove{D_T}\neq\emptyset$.\\
 Then, under assumption (C-0), we have
\begin{equation}
\label{lim3.Iinfty}
 \lim_{T'\nearrow T^*} \limsup_{\teta\in(0,T')} \lim_{\mu\to\infty}\lim_{\tau\to\infty}  \frac1{\tau}\ln |\, I_\infty(\tau,\mu,\teta,T')| = 0.
\end{equation}
\end{theor}
%

%
%

\medskip
 So we can detect if $\Sigma_T$ touches $D_T$ or not. Consider the function
$$
 F(T')= \limsup_{\teta\in(0,T')}\lim_{\mu\to\infty}\lim_{\tau\to\infty}\lim_{j\to\infty}
\frac1{\tau}\ln |\, I_j(\tau,\mu,\teta,T')|,\quad 0\le T'\le T,
$$
 defined in $[0,T^*[$. The problem is that $T^*$ is defined from $F$ itself.
 The following result gives a strategy to determine $T^*$.

\medskip
\begin{lem}
\label{l.c1}
 There exists $\delta>0$ such that if $T'<T^*$ then $T'<T'+|F(T')|/\delta<T^*$.
\end{lem}

\medskip
 In practice, if we have the a priori knowledge of $\delta$, we can determine $T^*$
 as the limit of the sequence
$$ (t_n)_{n\in\N}:\quad t_0=0;\quad t_{n+1} = t_n+ |F(t_n)|/\delta,\: n\ge 0. $$
 If $t_{n+1}> T$ for some $n$, this indicates that $\Sigma_T$ does not touch $D_T$.
 We write $T^*=T+0$ in this case.

\medskip
\begin{coro}
\label{c.1}
 From the knowledge of $\Lambda_{\gamma,v_0}$ we can compute :\\
 (1) the maximal time $T^*$ of a given curve $\Sigma$.\\
 (2) the connected components of $\ove{\Om_T}\sauf\ove{D_T}$ that touch
 $\{0\}\times\Gamma$.
\end{coro}

\medskip
 In particular, the connected components of $\ove{\Om_T}\sauf\ove{D_T}$ that
 touch $\{0\}\times\Gamma$ are completly characterized from $\Lambda_{\gamma,v_0}$.

 

\medskip 
 Let us make some remarks on Theorem \ref{t.1}, Corollary \ref{c.1}.

\medskip
\begin{remark}
 - We emphazise the fact that we don't make assumption on $D(t)$, $0\le t\le T$.\\
 - From the knowledge of $\Lambda_{\gamma,v_0}$,  we get the information not only
 on $\rmd(\Sigma_{T'},D(T'))$, when $T'<T^*$, but on $\rmd(y(\teta),D(\teta))$,
 $\teta\in (0,T')$, too.
\end{remark}

\medskip
\subsection{Literature review}
 Many articles solve a version of the Calder\'on inverse problem for the heat equation.
 The biggest part of them assume that the unknown coefficient  $\gamma$, or
 the unknown inclusion $D$, do not depend on time $t$: see for example
 \cite{IKE.ENC,IKE.REC}

 But they are very few results concerning the time dependent situation $D=D(t)$.
 The authors in \cite{ELA.UNI} proved uniqueness of $D_T$ under the assumption
 that the inclusion $D(t)$ is $x$-lipchitzian for all $t$. They used a proof
 by contradiction and is not constructive at all.
 A more recent paper for a similar question is \cite{KAW.UNI}.

 A reconstruction method by a dynamical probe method is performed
 in \cite{DAI.PRO}, but it works for the one dimensional spatial space only.
 In this case, an easier way is to put  as input the trace of plane
 waves $U_\tau(t,x)=e^{\tau^2 t+\tau x}$ at the spatial boundary, then
 to compute the solution $V_\tau(t,x)$ by using an ansatz and an energy estimate:
 see \cite{GAI.INV1}.
 This approach works in the $x$-multidimensional case if $D(t)$ is a convex set:
 see \cite{GAI.INV2}.
 
 With the aim of working the dynamical probe method in the $x$-multidimensional
 case with, the authors of \cite{ISA.REC} computed the reflected solution at points
 $(t,x)$ close to the lateral boundary of $D_T$. But this way is painful enough
 to assume that $D(t)=D$ does not depend on $t$. Moreover, since the data
 these authors consider are the traces on $\Gamma_T$ of Runge approximations of
 the fundamental solution $G_{s,y}$ of the heat equation, this method requires
 some regularity of $\partial D$ (more precisely, they assumed that $\partial D$
 is of class $C^{1,\alpha}$, for some $\alpha\in (0,1)$).

 In \cite{CRI.STA} the stability of the operator $\gamma\mapsto \Lambda_{\gamma,0}$
 is quantified, under the assumption that $D_T$ belongs to some class ${\mathcal K}$
 contained in $L^\infty((0,T)$; $W^{2,\infty}(\Om))$ $\cap$ $W^{1,\infty}((0,T)$;
 $W^{1,\infty}(\Om))$.

 We think that our approach, based on special fundamental solutions and
 inequalities\refq{dine.Iinfty}, would be able to extend the frameworks of
 many of these articles.
\noindent

\subsection{Outline of the paper}
 The paper is organized as follows.
 Section \ref{proofs} is devoted to the proofs and is devided in several subsections.
 In subsection \ref{s.BasEst}, we give basic estimates for integrals on $D_T$
 with weight functions as $(p_{\tau,\Sigma} \kappa)(t,x)$.
 In subsection \ref{s.SF}, we build the functions $U_\tau$, $U_\tau^*$
 by proving Lemma \ref{l.varphi} of section \ref{s.NS0}.
 In subsection \ref{ss.plca} we prove Lemma \ref{l.c1}.
 In subsection \ref{ss.pt1} we prove Theorem \ref{t.1}.
 In subsection \ref{ss.pt2} we prove Theorem \ref{t.2}.
 In subsection \ref{s.NS}, we build the needle sequences $U_j$, $U_{j,T'}$.
 
\section{Proofs}
\label{proofs}

\subsection*{Notations}
 We consider the usual Sobolev spaces $H^m(V)$, where $V\subset \R^k$ is an open set,
 and also $H^m((T_1,T_2);H^{m'}(V))=H^m(T_1,T_2;H^{m'}(V))$, for $m,m'\in\Z$.
 Here, $H^0$ is the well-know space $L^2$.
 
 The formal parabolic operator is $\calL_\gamma=\partial_t - \nabla_x\cdot(\gamma\nabla_x\cdot)$.

\subsection{Basic Estimates}
\label{s.BasEst}
 In this part we establish some basic estimates on $W_\tau$ etc.
%
%
%

\medskip
\begin{lem}
\label{l.BE}
a) Let $\calO\subset\R^3$ be an non empty open set and $\tau>0$.
 Then we have
\begin{eqnarray}
\label{est.py}
 \int_{\calO} |p_{\tau,y}(x)|^2 dx \le   C \frac1{\tau  \rmd(y,\calO)}
  e^{-2\tau \rmd(y,\calO)},  \quad \forall y\in\R^3\\
\label{max1.dpy}
 \int_{\calO} |\nabla_x p_{\tau,y}(x)|^2 dx \le  C\tau \left(1+\frac1{\tau \rmd(y,\calO)}\right)
   e^{-2\tau \rmd(y,\calO)}, \quad \forall  y\not\in\ove{\calO}.
\end{eqnarray}
b) 
 Let $\calO\subset\R^3$ is a non empty open. Let $y\not\in \ove\calO$ and
 $d>\rmd(y,\calO)$. Then there exists $C(\calO,y,d)>0$, such that
\begin{eqnarray}
\label{min1.dpS}
 \int_{\calO} |\nabla_x p_{\tau,y}(x)|^2 dx \ge C(\calO,y,d) \tau^2 e^{-2\tau d},
   \quad \forall \tau\ge 0.
\end{eqnarray}
\end{lem}

 Proof. a) We check easily:
\begin{eqnarray*}
 \int_{\calO} |p_{\tau,y}(x)|^2 dx &=& \int_{\calO} \frac1{(4\pi)^2 r^2} e^{-2\tau r} dx
  \le \int_{r>\rmd(y,\calO)} \frac1{4\pi} e^{-2\tau r} dr \\
 &\le & \frac1{\tau} \int_{s>\tau \rmd(y,\calO)} e^{-2s} ds
  = (2\tau \rmd(y,\calO))^{-1} e^{-2\tau \rmd(y,\calO)}.
\end{eqnarray*}
 Hence\refq{est.py}. Similarly we have
\begin{eqnarray*}
  \int_{\calO} |\nabla_x p_{\tau,y}(x)|^2 dx &=&
  \int_{\calO} \frac{(\tau r+1)^2}{(4\pi)^2 r^4} e^{-2\tau r} dx
  \le \int_{r>\rmd(y,\calO)} \frac{(\tau r+1)^2}{4\pi r^4} e^{-2\tau r} r^2 dr \\
  &\le & \frac{\tau}{4\pi} \int_{s>\tau \rmd(y,\calO)} \frac{(s+1)^2}{s^2} e^{-2s} ds\\
  &\le & C\tau \left(1+\frac1{\tau \rmd(y,\calO)}\right) e^{-2\tau \rmd(y,\calO)}.
\end{eqnarray*} 
b) Let us prove\refq{min1.dpS}. Observe that the open set $\calO_d:=\{x\in\calO;\; |x-y|<d\}$
 is non empty. Hence we get
\begin{eqnarray*}
  J(y) &:=& \int_{\calO} |\nabla_x p_{\tau,y}(x)|^2 dx \ge  \int_{\calO_d} \frac{(\tau r+1)^2}{(4\pi)^2 r^4}  e^{-2\tau d} dx
 \ge  C \tau^2 e^{-2\tau d},
\end{eqnarray*}
 with $C= \int_{\calO_d} \frac1{(4\pi r)^2} dx>0$. 
\begin{remark}
 It is clear that if $y\in\calO$, or if $\partial\calO$ is lipschitz and $y\in\partial\calO$,
 then
$$   \int_{\calO} |\nabla_x p_{\tau,y}(x)|^2 dx = +\infty .$$
\end{remark}
%
%

\medskip
\begin{lem}
\label{l.BS2}
 Assume that $\Sigma_T\cap \ove{D_T}=\emptyset$. Let $\teta\in [0,T]$.
 Let  $d> \rmd(y(\teta),D(\teta))$. Then, there exists $\mu_1>0$, $\tau_1>0$
 such that, for all $\mu>\mu_1$, $\tau>\tau_1$, we have
\begin{equation}
\label{max2.dpS}
 \int_{\Om_T} |\gamma-1|\,|\nabla_x p_{\tau,\Sigma}(x)|^2 dx\kappa(t)dt \le 
  \frac{C}{\mu} \big(1+\frac1{\tau \rmd(y(\teta),D(\teta))} \big)e^{-2\tau \rmd(y(\teta),D(\teta))},
\end{equation}
\begin{equation}
\label{min2.dpS}
 \int_{\Om_T} |\gamma-1|\,|\nabla_x p_{\tau,\Sigma}(x)|^2 dx\kappa(t)dt
  \ge   C\frac{\tau}{\mu} e^{-2\tau d},
\end{equation}
 for some $C>0$ which does not depend on $\mu$ or $\tau$.
\end{lem}
 Proof.

a) Thanks to\refq{max1.dpy} we have
\begin{eqnarray*}
 I &:=& \int_{\Om_T} |\gamma-1|\,|\nabla_x p_{\tau,\Sigma}(x)|^2 dx\kappa(t)dt \\
 &\le & C\tau \int_0^T   \Big(1+\frac1{\tau \rmd(y(t),D(t))} \Big) e^{-2\tau  \rmd(y(t),D(t))}
 \kappa(t) dt .
\end{eqnarray*}
 By assumption on $D_T$ and $\Sigma$, we have
%
$$
  \rmd(y(\teta),D(\teta)) + \delta |t-\teta| \ge \rmd(y(t),D(t))   \ge  \rmd(y(\teta),D(\teta))
 - \delta |t-\teta|,\quad t\in [0,T],
$$
%
 for some $\delta>0$ depending on $\Sigma_T$ and $D_T$ only.
 Fix $\mu_1>2\delta$ and let $\mu\ge\mu_1$.
 Put $\epsl=\frac{\rmd(y(\teta),D(\teta))}{\sqrt{\tau}}$. We have
\begin{eqnarray*}
 I &\le& C \tau  e^{-2\tau \rmd(y(\teta),D(\teta))}
  \Big(1+ \frac1{\tau (\rmd(y(\teta),D(\teta)) - \delta\epsl)} \Big)
   \int_{|t-\teta|<\epsl} e^{\tau(2\delta-\mu) |t-\teta|} dt \\
  && +  T e^{-2\tau \rmd(y(\teta),D(\teta))}
 e^{\tau\epsl(2\delta-\mu)}  \Big(1+ \frac1{\tau (\rmd(\Sigma_T,D_T)} \Big)  dt \\
 &\le& e^{-2\tau \rmd(y(\teta),D(\teta))} \Big\{ \frac{C}{\mu-2\delta} 
  \big(1+ \frac2{\tau \rmd(y(\teta),D(\teta))} \big) +  \calO(e^{-\delta'\sqrt{\tau}}) \Big\}.
\end{eqnarray*}
 for some $\delta'>0$. Estimate\refq{max2.dpS} is proved.
  
b) Put
$$
  J(t) :=\int_\Om |\gamma(t,x)-1|\, |\nabla_x p_{\tau,y(t)}(x)|^2 dx.
$$
 Thanks to assumption (C-D) we have
$$ J(\teta) \ge \delta \int_\Om |\nabla_x p_{\tau,y(\teta)}(x)|^2 dx,$$
 with $\delta=\liminf_{x\in \ove{\calO}} |\gamma(\teta,x)-1|>0$.
 Then, thanks to\refq{min1.dpS}, we have
$$ J(\teta)   \ge \delta'\tau^2 e^{-2\tau d},$$
 for some $\delta'>0$. Obviously, $J(\cdot)$ depends continuously on $t$, since
 $\gamma\in L^\infty(\Om_T)$.
 Hence if $|t-\teta|<\epsl$ and $\epsl>0$ is sufficiently small, we have
$$ J(t)  \ge \frac12\delta ' \tau^2 e^{-2\tau d}.$$
 We then have
\begin{eqnarray*}
 I & :=& \int_{\Om_T}J(t) \kappa(t)dt 
  \ge \int_{|t-\teta|<\epsl} J(t) \kappa(t) dt  \ge  C\tau^2 e^{-2\tau d}
  \int_{|t-\teta|<\epsl} \kappa(t)dt \\
  &\ge & C\frac{\tau}{\mu} e^{-2\tau d}(1-  e^{-\tau \mu\epsl}).
\end{eqnarray*}
 for some $C>0$. Taking $\mu_1>\frac1{\epsl}$, this proves\refq{min2.dpS}.

%
%
%

\medskip
\begin{lem}
 Let $\teta\in (0,T)$. Then there exist $\delta>0$, $\mu_1>0$ such that for all $\mu>\mu_1$
 there exist $C(\mu)\ge 0$, $\tau_1>0$ such that
\begin{eqnarray}
\label{est.U*t=0}
\nonumber
 \int_\Om  \kappa(0) |U|^2 dx &\le& C(\mu) e^{-2\tau \rmd(y(t),D(t))+\delta)},
  \quad t\in [0,T],\: \tau>\tau_1, \\
&& \quad U=U_\tau(0),U_\tau^*(0);
\end{eqnarray}
\begin{equation}
\label{est.U*t=T}
 \int_\Om  \kappa(T)|u_\tau^*(T)|^2 dx \le
  C(\mu) e^{-2\tau (\rmd(y(\teta),D(\teta))+\delta) }, \quad \tau>1;
\end{equation}
\begin{equation}
\label{est.Wt=0}
 \int_\Om \kappa(0) |W_\tau(0)|^2 dx \le  C(\mu) e^{-2\tau \rmd(y(t),D(t))+\delta) },
 \quad t\in [0,T],\: \mu>\mu_1,\: \tau>1;
\end{equation}
\begin{equation}
\label{est.WU*t=0}
\int_\Om  \kappa(0) |W_\tau(0) U^*(0)| dx \le C(\mu)  e^{-2\tau \rmd(y(t),D(t))+\delta) },
  \quad t\in [0,T],\: \tau>\tau_1.
\end{equation}
\end{lem}
 Proof.
a) Thanks to Lemma \ref{l.varphi},\refq{est.varphi}, to Remark \ref{r.p.U*},
 to\refq{est.py}, and reminding that $y(0)\not\in \ove\Om$, we have
\begin{eqnarray*}
 \int_\Om |U|^2 dx & \le & C \int_\Om |p_{\tau,y(0)}|^2 dx \le  C\tau^{-1} e^{-2\tau d(y(0),\Om)} \le C'.
\end{eqnarray*} 
 Since $\kappa(0)=e^{-\tau\mu\teta}$, we then have
\begin{eqnarray}
 \int_\Om  \kappa(0) |U|^2 dx \le  C'e^{-\tau\mu\teta}.
\end{eqnarray}
 Choosing $\mu_1> \sup_{(t,x)\in\Om_T} \frac{2\rmd(y(t),x)}{\teta}$, we obtain\refq{est.U*t=0}.\\
b) Similarly, choosing $\mu_1> \sup_{(t,x)\in\Om_T} \frac{2\rmd(y(t),x)}{(T-\teta)}$, we obtain\refq{est.U*t=T}.\\
c) Let us remind that $W_\tau(0,x)=v_0(x)-U_\tau(0,x)$.
 Thanks to Assumption (C-0) and to\refq{est.varphi},\refq{d.Utau},\refq{est.U*t=0},
 we have
$$ 
 \int_\Om \kappa(0) |W_\tau(0)|^2 dx \le Ce^{-\tau\mu\teta}(e^{l_0\tau}
 + \tau^{-1} e^{-2\tau d(y(0),\Om)})\le Ce^{-\teta(\mu-\mu_{T,\teta})\tau}.
$$
 Choosing $\mu_1> \mu_{T,\teta}+\sup_{(t,x)\in\Om_T} \frac{2\rmd(y(t),x)}{\teta}$
 we obtain\refq{est.Wt=0}.\\
d) Estimates\refq{est.U*t=0},\refq{est.Wt=0} and standart inequalities imply\refq{est.WU*t=0}.


\subsection{Construction of the special functions}
\label{s.SF}
 Let $R>0$ such that $\Om\subset \{x\, |\; |x|< R\}=B_R$.
 Put $z=x-y(t)$, $u_\tau(t,x)=f(t,z)$. 
 We have
$$ \partial_t u_\tau= \partial_t f + \dot z \nabla_z f= \partial_t f - \dot y \nabla_z f,$$
%
%
$$
 (\partial_t -\Delta_z)f  - \dot y(t) \nabla_z f + (\tau^2 + \tau\eta(t)) f = \delta_0(z).
$$
%
 Let $\hat f(t,k)= (2\pi)^{-3/2} \int_{z\in\R^3} f(t,z)e^{-ikz}dz$ be the Fourier transform
 of $f(t,\cdot)$.
 It satisfies
%
$$
 \partial_t \hat f + (k^2 - ik\dot y(t)+\tau^2 + \tau\eta(t)) \hat f  = 1.
$$
%
 Put $\rho(t)=\int_0^t \eta(s)ds$. We then have
\begin{equation}
\label{v.hf}
 \hat f(t,k) = (2\pi)^{-3/2} \!\!\! \int_{ \!\!-\infty}^t  \!\!\! \!\!\! \exp \left( k^2(s-t) - ik (y(s)-y(t))
  + \tau^2 (s-t)  + \tau(\rho(s)-\rho(t)) \right) ds,
\end{equation}
 and so
\begin{eqnarray}
\nonumber
 f(t,z) &=&  (2\pi)^{-3/2}\int_{-\infty}^t \left( \int_{k\in\R^3} \!\! \!\!e^{k^2(s-t)
  + ik(z+y(t)-y(s))} dk \right)
  e^{\tau^2 (s-t) + \tau(\rho(s)-\rho(t))} ds,\\
\nonumber
 &=& \frac1{8\pi^{3/2}} \int_0^{\infty} \left( \int_{k\in\R^3} \!\! \!\!e^{-k^2 s + ik(z+y(t)-y(t-s))} dk \right)
  e^{-\tau^2 s + \tau(\rho(t-s)-\rho(t))} ds,\\
\label{v.utau}
 u_\tau(t,x) &=& \frac1{8\pi^{3/2}} \int_0^{\infty}  e^{-\tau^2 s + \tau(\rho(t-s)-\rho(t))}  e^{\frac{-(x-y(t-s))^2}{4s}} s^{-3/2} ds.
\end{eqnarray}
 The integral in\refq{v.utau} is convergent if and only if $x\neq y(t)$. We set
$$ \varphi := \frac{u_\tau}{p_{\tau,\Sigma}}. $$
 Let us prove\refq{est.varphi}. We put
\begin{eqnarray*}
 g_\tau(s,t,x) &=& e^{\tau(\rho(t-s)-\rho(t))+ \frac1{4s}(y(t-s)-y(t))(2x-y(t-s)-y(t))},\\
 k_\tau(s,t,x) &=& \frac1{8\pi^{3/2}}   e^{-\tau^2 s}  e^{\frac{-(x-y(t))^2}{4s}} s^{-3/2},
\end{eqnarray*}
 and write
$$
 u_\tau(t,x) = \int_0^{\infty} k_\tau(s,t,x) g_\tau(s,t,x) ds.
$$
%
 Let us observe that
$$ p_{\tau,y(t)}(x)   = \int_0^{\infty} k_\tau(s,t,x)  ds,\quad x\in\R^3. $$ 
 In fact, the right-hand side in the above relation belongs to $L^2(\R^3)$, with a Fourier
 transform corresponding to\refq{v.hf} where $g$ is replaced by $1$.
 Moreover, it satisfies\refq{e.ptau}, since it is a time-independent solution of\refq{p.U}.
 Hence this right-hand side is $p_{\tau,y(t)}(x)$.

\medskip
 Since $\rho(\cdot)$ is lipschitzian, we can fix $K\in\R$ such that $\tau>|K|>\|\dot \rho\|_\infty$.
 Let us remind that $\dot\rho=\eta$ and $\eta$ satisfies\refq{est.eta}.
 Thus the above inequalities impose $\mu <\tau$, and, in fact,
$$ \mu <<\tau,\quad {\rm as}\quad \tau\to\infty .$$
 Then, if $K>0$, we can write, for $x\in \Om\subset B_R$, $x\neq y(t)$,
\begin{eqnarray*}
 k_\tau(s,t,x) g_\tau(s,t,x) &=&  k_{\tau-K}(s,t,x) e^{2\tau K s+K^2s} g_\tau(s,t,x) ds\\
  &\le & C(\mu,K,R) k_{\tau-K}(s,t,x).
\end{eqnarray*}
 Hence 
\begin{eqnarray*}
 u_\tau &\le& C(\mu,K,R) p_{\tau,\Sigma}\quad {\rm in}\quad \Om_T.
\end{eqnarray*}
 For $K<0$ we obtain
%
$$
 u_\tau(t,x) \ge p_{\tau-K,y(t)}(x) \ge C'(\mu,K,R) p_{\tau,y(t)}(x),
$$
%
 where $C'(\mu,K,\Om) >0$. Thus\refq{est.varphi} is proved.
 Let us prove\refq{est.dvarphi}. Let $x\in \Om\subset B_R$, $x\neq y(t)$,
 We use the following changes of variable
$$
 s=\frac{r}{2\tau}e^a,\quad \tilde s:=\frac{r}{2\tau}e^{-a}, \quad b=\cosh a\:
  (a>0), \quad \alpha=\tau r b',\quad b'=b-1.
$$
 We then have
\begin{eqnarray*}
 u_\tau(t,x) &=& \frac1{8\pi^{3/2}} \int_0^{\infty}   e^{-\tau^2 s-\frac{r^2}{4s}} s^{-3/2} g_\tau(s,t,x) ds \\
  &=& \frac1{8\pi^{3/2}} (\frac{r}{2\tau})^{-1/2} \int_{-\infty}^{\infty}   e^{-\tau r \cosh a} e^{-a/2} g_\tau(s,t,x) da\\
  &=& \frac1{8\pi^{3/2}} (\frac{2\tau}{r})^{1/2}  \int_0^{\infty}   e^{-\tau r \cosh a} (e^{-a/2} g_\tau(s,t,x)
   + e^{a/2} g_\tau(\tilde s,t,x))da \\
 &=& \frac1{8\pi^{3/2}} (\frac{2\tau}{r})^{1/2}  \int_0^{\infty}   e^{-\tau r b} \Big( \cosh (a/2) (g_\tau(s,t,x) + g_\tau(\tilde s,t,x)) -\\
 && \sinh (a/2) (g_\tau(s,t,x) - g_\tau(\tilde s,t,x)) \Big) da\\
 &\equiv & u^+ - u^-,
\end{eqnarray*}
 with
\begin{eqnarray*}
 u^+_\tau(t,x)  &=&  \frac1{8\pi^{3/2}} (\frac{2\tau}{r})^{1/2}  \int_0^{\infty}
  (g_\tau(s,t,x) + g_\tau(\tilde s,t,x)) e^{-\tau r b} \cosh (a/2)  da \\
  &=&  \frac1{8\pi^{3/2}} (\frac{2\tau}{r})^{1/2} e^{-\tau r} \int_0^{\infty}
  (g_\tau(s,t,x) + g_\tau(\tilde s,t,x))  e^{-\tau r b'}  \frac{\cosh (a/2)}{\sinh(a)} db' \\
  &=&  \frac{\sqrt2}{8\pi^{3/2}}  \frac1r e^{-\tau r} \int_0^{\infty} (g_\tau(s,t,x)
  + g_\tau(\tilde s,t,x))  \frac{\cosh (a/2)}{\sinh(a)} e^{-\alpha} \frac{d\alpha}{\sqrt{\tau r}}\\
 &\equiv& p_{\tau,y(t)} \varphi^+(t,x),\\
 u^-_\tau(t,x)  &=&  \frac1{8\pi^{3/2}} (\frac{2\tau}{r})^{1/2}  \int_0^{\infty}
  (g_\tau(s,t,x) - g_\tau(\tilde s,t,x)) e^{-\tau r b} \sinh (a/2)  da\\
  &=&  \frac{\sqrt2}{8\pi^{3/2}}  \frac1r e^{-\tau r} \int_0^{\infty} (g_\tau(s,t,x)
  - g_\tau(\tilde s,t,x)) \frac{\sinh (a/2)}{\sinh(a)} e^{-\alpha}  \frac{d\alpha}{\sqrt{\tau r}} \\
 &\equiv& p_{\tau,y(t)} \varphi^-(t,x).
%
%
\end{eqnarray*}
 Estimate\refq{est.dvarphi} in Lemma \ref{l.varphi} is a consequence of the relation $\varphi=\varphi^-+\varphi^+$
 and the folllowing Lemma.

\medskip
\begin{lem}
\label{l.vp+-}
 We have, for $r\le R$ and $\tau$ sufficiently large ($\tau>C(\mu,R)$),
\begin{eqnarray}
\label{est.vp-}
 |\varphi^-(t,x)| &\le &C(\mu,R)/\sqrt{\tau},\\
\label{est.vp+}
 |\varphi^+(t,x)| &\le & C(\mu,R),\\
\label{est.dvp-}
 |\nabla_x \varphi^\pm(t,x)| &\le &  C(\mu,R).
\end{eqnarray}
\end{lem}
 Proof.
 We observe that $0<\tilde s<1<s$ for $a>0$, and
\begin{eqnarray*}
  \cosh(a/2)&=& \sqrt{(b'+2)/2} = \sqrt{\frac{\alpha}{2\tau r}+1},\\
  \sinh(a/2)&=& \sqrt{b'/2} = \sqrt{\frac{\alpha}{2\tau r}},\\
\nonumber
  \tau (s-\tilde s) &=& r \sinh(a)  = r\sqrt{b^2-1} = r\sqrt{b'(b'+2)} =  \frac1{\sqrt{\tau}}\sqrt{\alpha}\sqrt{\frac{\alpha}{\tau}+2r},\\
  &\le&  \frac{\alpha}{\tau} + \sqrt{\frac{2r}{\tau}}\sqrt{\alpha},\\
 \tau(s+\tilde s) &=& r\cosh(a) = rb = r+rb' = r+ \frac{\alpha}{\tau}.
\end{eqnarray*}
 Hence, for $0<r\le R$, we have
\begin{eqnarray*}
 g_\tau(s,t,x) &=& {\rm exp} \Big\{-\tau (\rho(t-s)-\rho(t)) + \frac1{4s} (y(t-s)-y(t))(2r+y(t)-y(t-s)) \Big\} \\
 &\le&  {\rm exp} \Big(\tau s |\dot\rho|_\infty +  |\dot y|_\infty (r+ s|\dot y|_\infty)\Big)
 \le  e^{C(r+ \tau s)} \le  C e^{2Cr+C\frac{\alpha}{\tau}}
 \le  C e^{C\frac{\alpha}{\tau}},
\end{eqnarray*}
 where $C=C(\mu,R)$.
 The same estimate for $g_\tau(\tilde s,t,x)$ holds:
%
$$
 \max(g_\tau(s,t,x),g_\tau(\tilde s,t,x)) \le  C e^{C\frac{\alpha}{\tau}}.
$$
%
 We thus have
\begin{eqnarray*}
 |g_\tau(s,t,x) - g_\tau(\tilde s,t,x)| &\le &
 \Big(\tau|\rho(t-s)-\rho(t-\tilde s)| + \Big|\frac1{4s}(y(t-s)-y(t)) \\
 && \hspace{-3.7cm}
 \cdot(2x- y(t-s)-y(t)) -  \frac1{4\tilde s}(y(t-\tilde s)-y(t))(2x-y(t-\tilde s)-y(t)) \Big| \Big) \\
  &&  \hspace{-3.7cm}  \cdot \max(g_\tau(s,t,x),g_\tau(\tilde s,t,x)) \\
  &&  \hspace{-4cm} \le  \Big( |\dot\rho|_\infty \tau (s-\tilde s) +\frac14 |\dot y|_\infty
  ( 2r+ |\dot y|_\infty s)  + \frac14|\dot y|^2_\infty (s+\tilde s)\Big)   e^{C\frac{\alpha}{\tau}} \\
  &&  \hspace{-4cm} \le  C (\tau (s-\tilde s) + r+s+\tilde s)  e^{C\frac{\alpha}{\tau}} 
 \le  C(\frac{\alpha}{\tau} + \sqrt{\frac{2r}{\tau}}\sqrt{\alpha} + r 
 + \frac1{\tau}(r+ \frac{\alpha}{\tau}))  e^{C\frac{\alpha}{\tau}}\\
  && \hspace{-4cm}\le 
     C(\frac{\alpha}{\tau} + r) e^{C\frac{\alpha}{\tau}},\quad C=C(\mu,R).
\end{eqnarray*}
 Consequently, since
\begin{eqnarray*}
 \varphi^-(t,x) &=&  \frac1{\sqrt{2\pi}}  \int_0^{\infty} (g_\tau(s,t,x) - g_\tau(\tilde s,t,x))
  \frac{\sinh (a/2)}{\sinh(a)}  e^{-\alpha} \frac{d\alpha}{\sqrt{\tau r}}\\
 &=& \frac1{\sqrt{2\pi}}  \int_0^{\infty} (g_\tau(s,t,x) - g_\tau(\tilde s,t,x))
 \frac1{2\cosh(a/2)} e^{-\alpha}  \frac{d\alpha}{\sqrt{\tau r}},
\end{eqnarray*}
 we then have, for $\tau>C=C(\mu,R)$,
\begin{eqnarray*}
 |\varphi^-(t,x)| &\le & C \int_0^{\infty} (\frac{\alpha}{\tau} + r) e^{(\frac{C}{\tau}-1)\alpha}
  \frac{d\alpha}{\sqrt{\tau r} \sqrt{\frac{\alpha}{2\tau r}+1}} \\
 &\le &  C \int_0^{\infty} (\frac{\sqrt{\alpha}}{\tau} + \sqrt{\frac{r}{\tau}})
  e^{(\frac{C}{\tau}-1)\alpha}  d\alpha\\
 &\le &  C'  (\frac1{\tau}+\sqrt{\frac{r}{\tau}}).
\end{eqnarray*}
 Thus\refq{est.vp-} is proved.
 Moreover, since
\begin{eqnarray}
\nonumber
 \nabla_x g_\tau(s,t,x) &=& \frac{y(t-s)-y(t)}{2s} g_\tau(s,t,x),\\
\label{est.dg}
 |\nabla_x g_\tau(s,t,x)| &\le & \frac12 |\dot y|_\infty  |g_\tau(s,t,x)| \le   Ce^{C\frac{\alpha}{\tau}},
\end{eqnarray}
 we then have
\begin{eqnarray*}
  |\nabla_r \varphi^-(t,x)| &\le & C  \int_0^{\infty} |g_\tau(s,t,x)+ g_\tau(\tilde s,t,x)|
   \frac1{2\cosh(a/2)} e^{-\alpha} \frac{d\alpha}{\sqrt{\tau r}}\\
 &\le &  C \int_0^{\infty} e^{(\frac{C}{\tau}-1)\alpha}\frac1{2\sqrt{\frac\alpha{2\tau r}+1}} \frac{d\alpha}{\sqrt{\tau r}} \\
 &\le &  C \int_0^{\infty} e^{(\frac{C}{\tau}-1)\alpha} \frac{d\alpha}{\sqrt{\alpha}} 
 \le  C'.
\end{eqnarray*}
 Thus the estimate of $\nabla_r \varphi^-(t,x)$ in\refq{est.dvp-} is proved.
 Let us prove\refq{est.vp+}. Since
$$
  \frac{\cosh (a/2)}{\sinh(a)} =  \frac1{2\sinh(a/2)} =  \frac1{\sqrt{2b'}}
  = \frac{\sqrt{\tau r}}{\sqrt{2\alpha}},
$$
 we then have
\begin{eqnarray*}
 \varphi^+(t,x) &=&  \frac1{\sqrt{2\pi}}  \int_0^{\infty} (g_\tau(s,t,x) + g_\tau(\tilde s,t,x))
   \frac{\cosh (a/2)}{\sinh(a)}  e^{-\alpha}  \frac{d\alpha}{\sqrt{\tau r}},
\end{eqnarray*}
 and so
\begin{eqnarray*}
 |\varphi^+(t,x)| &\le &  C |\int_0^{\infty} e^{C\frac{\alpha}{\tau}-\alpha}
  \frac1{\sqrt{2\alpha}}  d\alpha \le C'.
\end{eqnarray*}
 This proves\refq{est.vp+}. Similarly, thanks to\refq{est.dg}, we have
\begin{eqnarray*}
 |\nabla_x \varphi^+(t,x)|  &=&  \frac1{\sqrt{2\pi}} \left| \int_0^{\infty} (\nabla_x g_\tau(s,t,x)
  + \nabla_x g_\tau(\tilde s,t,x))  \frac1{\sqrt{2\alpha}} e^{-\alpha} d\alpha \right|,\\
  &\le&  C \int_0^{\infty} e^{C\frac{\alpha}{\tau}-\alpha} \frac1{\sqrt{2\alpha}}  d\alpha
  \le C',
\end{eqnarray*}
 which proves\refq{est.dvp-} for $\nabla_x\varphi^+$. So Lemma \ref{l.vp+-} is proved.

\medskip
 By observing that $\varphi=\varphi^-+\varphi^+$, Lemma \ref{l.varphi} is then proved.

\subsection{Properties of reflected waves}
 Putting
\begin{eqnarray}
\nonumber
 \Psi_\tau &=& (\gamma-1)\nabla_x V_\tau  + \nabla_x W_\tau =
  \gamma\nabla_x V_\tau-\nabla_x U_\tau,\\
%
\label{r.Psijujwj}
\Psi_j &=& (\gamma-1)\nabla_x V_j  + \nabla_x W_j = \gamma\nabla_x V_j-\nabla_x U_j ,
\end{eqnarray}
 we have
\begin{eqnarray}
\label{r.divPsitau}
 \dive_x \Psi_\tau &=& \dive_x(\gamma(t,\cdot)\nabla_x V_\tau) -\Delta_x U_\tau
 = \partial_t W_\tau \quad{\rm in}\quad \Om_T,\\
\label{r.divPsij}
 \dive_x \Psi_j &=& \dive_x(\gamma(t,\cdot)\nabla_x V_j) -\Delta_x U_j  
 = \partial_t W_j \quad{\rm in}\quad \Om_T.
\end{eqnarray}
 If $\Sigma_T\cap \ove{D_T}=\emptyset$ then, in view of\refq{r.Psijujwj},\refq{lim.UjUj*},
 \refq{lim.Wj}, we have
%
%
$$
 \Psi_j\stackrel {j\to\infty}\to \Psi_\tau \quad {\rm in} \quad C([0,T],(L^2(\Om))^3).
$$
%

\subsection{Proof of Lemma \ref{l.0}}
 Let us set $T'=T$ which does not restrict the proof. We denote $U_j^*$ instead of $ U_{j,T}^*$.
 Put
\begin{eqnarray*}
 J_j(t) &: =& \int_{\partial\Om}(\partial_\nu W_j)(t,x) \, U_j^*(t,x) d\sigma(x),\\
 F_j(t) &:=&  \int_{\Om} (\gamma-1)\nabla_x V_j \nabla_x U_j^* dx.
\end{eqnarray*}
 Thanks to\refq{r.divPsij} and integration by parts, we have
\begin{eqnarray*}
 J_j(t) &=& \int_{\partial\Om} (\nu\cdot \Psi_j) U_j^* d\sigma(x)=
  \int_{\Om} \Psi_j \nabla_x U_j^* dx + \int_{\Om} (\dive_x \Psi_j) U_j^* dx \\
 &=& F_j(t) +  \int_\Om (\nabla_x W_j \nabla_x U_j^* + \partial_t W_j U_j^*) dx.
\end{eqnarray*}
 Thanks to\refq{bc.Wj},\refq{e.Uj*}, we compute
\begin{eqnarray*}
 \int_\Om \nabla_x W_j \nabla_x U_j^* dx = - \int_\Om W_j \Delta_x U_j^* dx
 =  \int_\Om W_j (\partial_t U_j^*+ \calL_1^* U_j^*)  dx\\
 = \int_\Om W_j (\partial_t U_j^* + \frac{\dot\kappa}{\kappa} U_j^*)  dx
 = \int_\Om W_j  \kappa^{-1} \partial_t (\kappa U_j^*)  dx .
\end{eqnarray*}
 Thus
\begin{eqnarray}
\label{r.JjFj}
 J_j(t) 
 &=& F_j(t) + \int_{\Om} \kappa^{-1}  \partial_t (W_j U_j^* \kappa)  dx.
\end{eqnarray}
 Thanks to\refq{e.Uj*},\refq{r.JjFj}, an integration by parts according to $t$ brings
\begin{equation}
\label{max1.Ij}
 I_j  \equiv  \int_0^T  J_j(t) \kappa(t) dt =  \int_0^T F_j(t) \kappa(t)  dt + \int_\Om [W_j  U_j^* \kappa]^T_0 dx,
\end{equation}

 Thanks to\refq{ci.Uj*},\refq{lim.Wj},\refq{lim.Vj},\refq{max1.Ij}, and to Remark \ref{rem.Uj*},
 if $\Sigma_T\cap \ove{D_T}=\emptyset$, then we have:
\begin{eqnarray*}
 \int_\Om [W_j  U_j^* \kappa]^T_0 dx &\stackrel{j\to+\infty}\to& \int_\Om [W_\tau  U_\tau^* \kappa]^T_{0} dx,\\
  \int_0^T F_j(t) \kappa(t)  dt &\stackrel{j\to+\infty}\to & \tilde I_\tau \equiv \int_{\Om_T} (\gamma-1)\nabla_x V_j \nabla_x U_j^*\kappa(t)  dt ,\\
 I_j  &\stackrel{j\to+\infty}\to&  I_\infty(\tau,\mu,\teta,T) \equiv \tilde I_\tau +
   \int_\Om [W_\tau  U_\tau^* \kappa]^T_{0} dx,
\end{eqnarray*}
 which proves\refq{lim1.Ij} in (1) of Lemma \ref{l.0}.

\subsection{Proof of Theorem \ref{t.1}}
\label{ss.pt1}
 For the sake of simplicity, we assume that $\gamma\ge 1$.
 The case $\gamma\le 1$ is similar.
 We put
$$ F_\tau := \int_\Om (\gamma-1)\nabla_x V_\tau \nabla_x U_\tau^* dx. $$

\subsubsection{Lower Bound for $I_\infty$}
\label{sss.LB}
 Let us give another expression of
%
$$ \tilde I_\tau = \int_0^T F_\tau(t,x) \kappa(t) dt. $$
%
 We remind that
\begin{equation}
\label{r.U*U}
 U_\tau^* =    \rho e^{-2\tau^2 t} U_\tau,\quad \rho := \frac{\varphi^*}{\varphi},  \end{equation}
 where $\varphi$, $\varphi^*$ are characterized in Lemma \ref{l.varphi},
 and we put
\begin{eqnarray*}
 r_1  &:=& \frac{\gamma-1}{\gamma} \Psi_\tau (\frac{\nabla_x \rho}{\rho}) U_\tau^*,\\
 r_2 &:=& \frac{\gamma-1}{\gamma} \Big \{ (\nabla_x \varphi^*\varphi +  \varphi^*\nabla_x\varphi )
   p_{\tau,\Sigma} \nabla_x p_{\tau,\Sigma}  + \nabla_x \varphi^* \nabla_x\varphi  |p_{\tau,\Sigma}|^2 \Big\}.
\end{eqnarray*}
   
 We firstly have
\begin{eqnarray*}
 (\gamma-1)\nabla_x V_\tau \nabla_x U_\tau^* &\stackrel{(1)}=&  \frac{\gamma-1}{\gamma} \Psi_\tau
  \nabla_x U_\tau^* +  \frac{\gamma-1}{\gamma} \nabla_x U_\tau \nabla_x U_\tau^*\\
 && \hspace{-2cm} \stackrel{(2)}=  r_1 + \frac{\gamma-1}{\gamma} \Psi_\tau\rho \nabla_x U_\tau e^{-2\tau^2 t}  
 + \frac{\gamma-1}{\gamma} \nabla_x U_\tau \nabla_x U_\tau^*  \\
 && \hspace{-2cm} \stackrel{(3)}=  r_1 + \frac1{\gamma} \Psi_\tau^2 \rho  e^{-2\tau^2 t}
 -  \Psi_\tau (\nabla_x W_\tau) \rho e^{-2\tau^2 t}
 + \frac{\gamma-1}{\gamma} \nabla_x U_\tau \nabla_x U_\tau^* \\
 && \hspace{-2cm} =  r_1 + \frac1{\gamma} \Psi_\tau^2 \rho  e^{-2\tau^2 t}
 -  \Psi_\tau (\nabla_x W_\tau) \rho e^{-2\tau^2 t}
 +\frac{\gamma-1}{\gamma} |\nabla_x p_{\tau,\Sigma}|^2 \varphi^*\varphi + r_2.
\end{eqnarray*}
 Explanations - (1): write $\nabla_x V_\tau = \frac1{\gamma}(\Psi +\nabla_x U_\tau)$. 
 (2): use\refq{r.U*U}. (3): write $(\gamma-1)\nabla_x U_\tau = \Psi_\tau -\gamma \nabla_x W_\tau$
  in the second term.

 By integration by parts, thanks to\refq{bc.Wtau} and\refq{r.divPsitau}, we have
\begin{eqnarray}
\label{r.nWPsi}
 - \int_\Om \Psi_\tau \nabla_x W_\tau \rho dx &=&  \int_\Om \partial_t W_\tau W_\tau \rho dx
 + r_3,\\
\nonumber
 r_3 &:=&  \int_\Om \Psi_\tau W_\tau \nabla_x \rho dx.
\end{eqnarray}
 We thus have
\begin{eqnarray*}
 F_\tau &=&   \int_\Om  \frac{\gamma-1}{\gamma} |\nabla_x p_{\tau,\Sigma}|^2 \varphi^*\varphi  e^{-2\tau^2 t}dx
   +  \int_\Om \frac1{\gamma} \Psi_\tau^2 \rho e^{-2\tau^2 t} dx \\
&&  +  \int_\Om (\partial_t W_\tau) W_\tau \rho  e^{-2\tau t^2} dx + \int_\Om ( \frac{\gamma-1}{\gamma} (r_1+r_2) + r_3)dx.
\end{eqnarray*}
 Writing
\begin{eqnarray}
\label{dtWW}
 \int_{\Om_T} (\partial_t W_\tau) W_\tau \rho  e^{-2\tau t^2}  \kappa dt dx &=&
  \int_{\Om_T}  (\tau^2 \rho- \frac12 \dot\rho) |W_\tau|^2  e^{-2\tau t^2} \kappa dt dx\\
\nonumber && + \frac12\int_\Om [W_\tau^2 e^{-2\tau t^2} \kappa]^T_0,
\end{eqnarray}
 we then have
\begin{eqnarray}
\nonumber
 & \hspace{-5.5cm}   I_\infty(\tau,\mu,\teta,T)  =  \tilde I_\tau + \int_\Om [W_\tau  U_\tau^* \kappa]^T_0 dx \\
\label{vmin.Iinfty}
 & \hspace{-2.2cm}  = \int_{\Om_T}  \frac{\gamma-1}{\gamma} |\nabla_x p_{\tau,\Sigma}|^2 \varphi^*\varphi  \kappa dt  dx + \int_{\Om_T} \frac1{\gamma} \Psi_\tau^2
 \rho e^{-2\tau^2 t}\kappa dt dx \\
\nonumber
 & \hspace{-0.2cm} + \int_{\Om_T}  (\tau^2 \rho- \frac12 \dot\rho)
 |W_\tau|^2  e^{-2\tau t^2} \kappa dt dx + \frac12 \int_\Om  (W_\tau^2 e^{-2\tau t^2})
  \kappa|_{t=T}dx +R,
\end{eqnarray}
 where we put
\begin{eqnarray*}
 R &:=& R_1+R_2+R_3 + R_4 + R_5,\\
 R_j &:=&   \int_{\Om_T}  r_j \kappa dt dx, \quad j=1,2,3,\\
 R_4 &:=& \int_\Om [W_\tau  U_\tau^* \kappa]^T_0 dx, \\
 R_5 &:=&  -  \frac12 \int_\Om (W_\tau^2 \kappa)|_{t=0} dx,
\end{eqnarray*}

\medskip
 We prove in the next Lemma that $R$ is neglictable.

\medskip
\begin{lem}
\label{maj.Rj}
 There exist $\delta>0$, $\mu_1>0$, so that for all $\mu>\mu_1$, there exists $\tau_1>0$
 so that for all $\tau>\tau_1$, the following estimates hold.
\begin{eqnarray}
\label{est.R1}
 |R_1| &\le&  \tau^{-1} \int_{\Om_T} \Psi_\tau^2 \rho e^{-2\tau^2 t}\kappa dt dx + 
   C\tau^{-1} \int_{\Om_T} |\gamma-1| \; |\nabla_x p_{\tau,\Sigma}|^2 \kappa dt dx;\\
\label{est.R2}
 |R_2| &\le& C\tau^{-1} \int_{\Om_T} |\gamma-1| \; |\nabla_x p_{\tau,\Sigma}|^2 \kappa dt dx;\\
\label{est.R3}
 |R_3| &\le&   \int_{\Om_T} \frac1{2\gamma} \Psi_\tau^2 \rho e^{-2\tau^2 t}\kappa dt dx + 
 C \int_{\Om_T} |W_\tau|^2 e^{-2\tau^2 t}\kappa dt dx;\\ 
\label{est.R4}
 |R_4| &\le&  \frac14 \int_\Om |W_\tau|^2|_{t=T}  e^{-2\tau T}\kappa(T) dx +
   C e^{-2\tau(y(t),D(t))+\delta) },\\ 
\label{est.R5}
 |R_5| &\le&  C e^{-2\tau(y(t),D(t))+\delta) },\quad t\in [0,T].
\end{eqnarray}
\end{lem} 
Proof.
a) Since $p_{\tau,\Sigma}\le \tau |\nabla_x p_{\tau,\Sigma}|$ and thanks to Remark \ref{r.p.U*},
 we obtain\refq{est.R1} and\refq{est.R2}.\\
b) Estimates\refq{est.R3} comes from\refq{est.varphi},\refq{est.dvarphi}, from
 Remark \ref{r.p.U*} and standart inequalities.\\
c) Thanks to Lemma \ref{l.varphi},\refq{est.varphi}, to Remark \ref{r.p.U*}, to\refq{est.py},
 we have
$$
 |R_4| \le  \frac14 \int_\Om |W_\tau|^2|_{t=T}  e^{-2\tau T}\kappa(T)) dx +
  C\int_\Om |p_{\tau,y(T)}|^2 \kappa(T) dx + \int_\Om |W_\tau U_\tau^*|_{t=0} \kappa(0) dx
$$
d) Estimate\refq{est.R5} is\refq{est.Wt=0}. Then\refq{est.R4} comes from\refq{est.U*t=0}
,\refq{est.WU*t=0}. This proves the Lemma.

\medskip
 From Lemma \ref{maj.Rj},\refq{vmin.Iinfty},\refq{est.varphi}, Remark \ref{r.p.U*}, and since
$$ \tau^2 \rho- \frac12 \dot\rho + C\ge C' \tau^2, \quad C'>0, $$
 for $\tau$ sufficiently large, we obtain
\begin{eqnarray}
\nonumber
  I_\infty(\tau,\mu,\teta,T)
 &\ge & \frac1C \Big\{ \int_{\Om_T}  \frac{\gamma-1}{\gamma} |\nabla_x p_{\tau,\Sigma}|^2  \kappa dt  dx 
   + \int_{\Om_T} \Psi_\tau^2 \rho e^{-2\tau^2 t}\kappa dt dx \\
\nonumber &&  + \int_{\Om_T} \tau^2 |W_\tau|^2  e^{-2\tau t^2} \kappa dt dx
  +  \int_\Om  |W_\tau|^2|_{t=T}  e^{-2\tau T^2}\kappa(T)dx \Big\},
\end{eqnarray}
 for some $C>1$, which implies the lower bound for $|I_\infty|$ in\refq{dine.Iinfty}.

\subsubsection{Upper Bound for $I_\infty$}
 In the same way, putting
$$  r_6 := (\gamma-1) \nabla_x W_\tau  (\nabla_x \rho) U_\tau e^{-2\tau^2 t}, $$
 we have,
\begin{eqnarray*}
 (\gamma-1)\nabla_x V_\tau \nabla_x U_\tau^* & \stackrel{(1)}=&
  (\gamma-1) \nabla_x W_\tau  \nabla_x U_\tau^* + (\gamma-1) \nabla_x U_\tau \nabla_x U_\tau^*  \hspace{2cm}\\
 &&  \hspace{-3cm} \stackrel{(2)}= (\gamma-1) \rho \nabla_x W_\tau  \nabla_x U_\tau e^{-2\tau^2 t}
  + (\gamma-1) \nabla_x U_\tau \nabla_x U_\tau^* +r_6 \\
 &&  \hspace{-3cm} = \gamma \rho \nabla_x W_\tau  \nabla_x U_\tau e^{-2\tau^2 t}
  -\rho \nabla_x W_\tau  \nabla_x U_\tau e^{-2\tau^2 t} + (\gamma-1) \nabla_x U_\tau \nabla_x U_\tau^* +r_6 \\
 && \hspace{-3cm} \stackrel{(3)}= -\gamma \rho \nabla_x |W_\tau|^2 e^{-2\tau^2 t} + \rho \nabla_x W_\tau \Psi_\tau e^{-2\tau^2 t}
  + (\gamma-1) \nabla_x U_\tau \nabla_x U_\tau^* +r_6.
\end{eqnarray*}
 Explanations - (1): write $V_\tau = U_\tau+W_\tau$. 
 (2): use\refq{r.U*U}. (3): write $U_\tau = V_\tau-W_\tau$ in the firs two terms.

 Thanks to\refq{r.nWPsi}, we have
\begin{eqnarray*}
 F_\tau &=&  -\int_\Om \gamma \rho \nabla_x |W_\tau|^2 e^{-2\tau^2 t} dx
   - \int_\Om \partial_t W_\tau W_\tau \rho  dx - r_3\\
 &&  +  \int_\Om (\gamma-1) \nabla_x U_\tau \nabla_x U_\tau^* dx +  \int_\Om r_6 dx.
\end{eqnarray*}
 Thanks to\refq{r.nWPsi},\refq{dtWW}, we obtain
\begin{eqnarray*}
 \tilde I_\tau  &=& - \int_{\Om_T} \gamma \rho |\nabla_x W_\tau|^2 e^{-2\tau^2 t} \kappa dt dx 
  - \int_{\Om_T}  (\tau^2 \rho- \frac12 \dot\rho) |W_\tau|^2  e^{-2\tau t^2} \kappa dt dx \\
 &&   - \frac12 \int_\Om (W_\tau^2 e^{-2\tau t^2} \kappa)_{t=T}dx
   +  \int_{\Om_T} (\gamma-1) \nabla_x U_\tau \nabla_x U_\tau^*\kappa dt dx \\
 && -R_3 -R_5+ R_6,
\end{eqnarray*}
 where we put
$$ R_6 := \int_{\Om_T} r_6 \kappa(t) dt dx. $$
 Similarly to the above section \S \ref{sss.LB}, for $\tau$ sufficiently large, we obtain
\begin{eqnarray}
\nonumber
  I_\infty(\tau,\mu,\teta,T)
 &\le & \frac1C \Big\{ \int_{\Om_T}  (\gamma-1) |\nabla_x p_{\tau,\Sigma}|^2  \kappa dt  dx 
   - \int_{\Om_T} |\nabla_x W_\tau|^2 e^{-2\tau^2 t} \kappa dt dx \\
\nonumber &&  - \int_{\Om_T} \tau^2 |W_\tau|^2  e^{-2\tau t^2} \kappa dt dx
  - \int_\Om  |W_\tau|^2|_{t=T}  e^{-2\tau T^2}\kappa(T)dx \Big\},
\end{eqnarray}
 for some $C>1$, which implies the upper bound for $|I_\infty|$ in\refq{dine.Iinfty}.
 
\medskip
 Thus Theorem \ref{t.1} is proved.

\medskip
\subsection{Proof of Theorem \ref{t.2}}
\label{ss.pt2}


a) Assume that $T<T^*$. Let $\teta\in (0,T)$ and $\mu >\mu_1$, where $\mu_1$
 is defined in Theorem \ref{t.1}.
 Thanks to\refq{min2.dpS} and to Theorem \ref{t.1}, we obtain\refq{lim1.Iinfty}.
 Since in this formula the right-hand side does not depend on $\mu$, we let
 $\mu$ tend to infinity. We then obtain\refq{lim2.Iinfty}.

b) Assume that  $\Sigma_T\cap \ove{D_T}\neq\emptyset$, that is, $T\ge T^*$.
 Let $0<T'<T^*$. Thus,\refq{lim2.Iinfty} holds with $T$ replaced by $T'$.
 Taking the limit as $T'\nearrow T^*$, and since $\rmd(\Sigma_{T'},D_{T'})$ tends
 to $\rmd(\Sigma_{T^*},D_{T^*}) =0$ when, we obtain\refq{lim3.Iinfty}.

\subsection{Proof of Lemma \ref{l.c1}}
\label{ss.plca}
 First remind that $F$ is defined at least for $T'$ sufficiently close to 0.
 Observe that the function $G$:
$$ G:\quad [0,T]\supset T'\mapsto -2\rmd(\Sigma_{T'},D_{T'}) $$
 is lipschitzian, non-decreasing, and $G(0)<0$. Thus, there exists $\delta>0$
 such that
$$ G(T'') \le G(T') + 2\delta (T''-T'),\quad 0\le T'\le T''\le T. $$
 Hence, if $G(T')<0$, then $G(T'+|G(T')|/\delta)<0$.
 The result is proved since , $G=F$ on $[0,T^*[$.

\subsection{Existence of Needle sequences}
\label{s.NS}
 Denote by $\Sigma^C$ the open set $((-1,T+1)\times\Om' )\sauf \Sigma_{[-1,T+1]}$,
 with $\Sigma_{[t,t']}=\{(s,y(s));\, t\le s\le t'\}$.
 Let us prove the existence of $\{u_j\}_j$. Let $G(t,x)$ be the solution in $L^2((0,T+1);H^1_0(\Om'))
 \cap H^1((0,T+1);L^2(\Om'))\subset C([-1,T+1];L^2(\Om'))$ of
\begin{eqnarray*}
  (\calL_1 +\tau \eta)G &=&0 \quad {in}\quad \Om'_{T+1},\\
 G|_{t=0} &=& U_\tau(0,\cdot) \quad {in}\quad \Om',
\end{eqnarray*}
 and $\partial_t G\in L^2((0,T+1);H^{-1}(\Om'))$. Since $U_\tau\in C([-1,T+1];L^2(\Om'))$ and
 $\partial_t U_\tau\in L^\infty((-1,T+1);H^{-1}(\Om'))$, then the function $H=U_\tau-G$ extended by 0
 for $t<0$ belongs to $C([-1,T+1];L^2(\Om'))$ with $\partial_t H\in L^2((-1,T+1);H^{-1}(\Om'))$,
 and satisfies
\begin{eqnarray*}
 (\calL_1 +\tau \eta) H \,(t,x)&=& \chi_{(0,+\infty)}(t)\otimes\delta(x-y(t)) \quad {in}\quad (-1,T+1)\times\Om',\\
 H|_{t\le 0} &=& 0 \quad {in}\quad \Om'.
\end{eqnarray*}
 Let $U$ be a simply connected open set in $\R^4$ with lipschitz boundary, and such that $\ove U\subset \Sigma^C$.
 Observe that $\Sigma^C \sauf \ove U$ is connected. For all open set $V$ such that $\ove U\subset V\subset\ove V
 \subset \Sigma^C$, we have $H|_V\in H^{1,2}(V)$ , $H|_V(t,\cdot)=0$ for $-1<t\le 0$, and $(\calL_1 +\tau \eta) H =0$
 in $V$.
 Similarly to \cite[Theorem 2.3]{DAI.PRO},  there exists a sequence $H_j\in (-1,T+1)\times\Om'$ such that
\begin{eqnarray}
\label{e.Hj}
 (\calL_1 +\tau \eta) H_j \,(t,x)&=& 0 \quad {\rm in}\quad (-1,T+1)\times\Om',\\
\nonumber
 H_j|_{t\le 0} &=& 0 \quad {\rm in}\quad \Om',
\end{eqnarray}
 and $H_j$ converges to $H$ in $L^2(U)$ (with the strong norm).
 We then have the existence of another sequence $H'_j\in H^{1,2}_{loc}(\Sigma^C)$ satisfying
 the homogeneous relation\refq{e.Hj} and such that $H'_j$ converges to $H$ in $H^{1,2}_{loc}(\Sigma^C)$,
 as noted in \cite[\S3]{DAI.PRO}.
 Finally, the sequence $U_j=H_j+G$ satisfies Equation\refq{e.Hj} in $\Om'_T$, with
$U_j|_{t=0} =U_\tau(0)$ in $\Om'$, and $U_j\to U_\tau$ in $H^{1,0}_{loc}(\Om_T\sauf\Sigma_T)\cap H^{0,1}_{loc}(\Om_T\sauf\Sigma_T)$.

\medskip
\begin{remark}
 We used the fact that $U_\tau(0)\in L^2(\Om')$ only. Hence, the proof for the existence
 of $U_j^*$ in the dual case, where the triplet $(\calL_1,t,t=0)$ is replaced by
 $(\calL_1^*,T-t,t=T)$, works similarly.
\end{remark} 

\end{document}